\documentclass{article}
\usepackage{amsmath}
\usepackage{amssymb}
\usepackage{amsthm}
\usepackage{float}
\newtheorem{Theorem}{Theorem}
\newtheorem{Definition}{Definition}

\newtheorem{Lemma}{Lemma}

\usepackage{paralist}
\usepackage{graphicx}

\newcommand{\Includegraphics}[2]%
{\ifFigs\includegraphics[      width=#1]{#2}%
	\else  \includegraphics[draft,width=#1]{#2}\fi}

\begin{document}
	
	\title{Dynamics of 1D discontinuous maps \\ with multiple partitions 
		and linear functions  \\ having the same fixed point.  \\ 
		An application to financial market modeling}
	\author{\textit{Laura Gardini}$^{1,2},$ 
		\textit{ Davide Radi}$^{2,3}$, 
		\textit{Noemi Schmitt}$^{4}$, \\ 
		\textit{Iryna Sushko}$^{3,5}$, 
		\textit{Frank Westerhoff}$^{4}$\\ 
		$^{1}${\small Dept of Economics, Society and Politics, University of Urbino
			Carlo Bo, Italy}\\ 
		$^{2}${\small Dept of Finance, V\v{S}B - Technical University of Ostrava,
			Ostrava, Czech Republic}\\ 
		$^{3}${\small Dept of Mathematics for Economic, Financial and Actuarial
			Sciences,} \\ {\small Catholic University of Milan, Italy} \\ 
		$^{4}${\small Dept of Economics, University of Bamberg, Germany}\\ 
		$^{5}${\small Inst. of Mathematics, NAS of Ukraine, Kyiv, Ukraine}}
	
	\date{}
	\maketitle
	
	\textbf{Abstract}
	\smallskip

Piecewise smooth systems are intensively studied today in many application
areas, such as economics, finance, engineering, biology, and ecology. In this
work, we consider a class of one-dimensional piecewise linear discontinuous
maps with a finite number of partitions and functions sharing the same real
fixed point. We show that the dynamics of this class of maps can be analyzed
using the well-known piecewise linear circle map and prove that their bounded
behavior, unrelated to the fixed point, may consist of either nonhyperbolic
cycles or quasiperiodic orbits densely filling certain segments, with possible
coexistence. A corresponding model describing the price dynamics of a
financial market serves as an illustrative example. While simulated model
dynamics may be mistaken for chaotic behavior, our results demonstrate that
they are quasiperiodic.

\smallskip


\section{Introduction}

Since the seminal work by Huang and Day \cite{HD-93} on bull and bear market
dynamics, numerous models have been published that use one-dimensional (1D)
piecewise linear (PWL) or piecewise smooth (PWS) maps to explain the behavior
of financial markets (see, e.g., \cite{H10,TWG-10,TWG-15,Brianzoni,Lu}). In these models, different types of market
participants interact. The key actors shaping market dynamics, besides the
market maker---who adjusts prices in response to excess demand---are
fundamentalists, who buy or sell assets when the market is undervalued or
overvalued, and chartists, who buy or sell assets when they perceive a bull or
bear market. As has been shown, even simple low-dimensional deterministic
models can capture the intricate price behavior of financial markets,
exhibiting transitions from fixed-point dynamics to chaotic behavior under
parameter changes. 

The goal of our paper is twofold. First, we contribute to the line of research
that studies the behavior of financial markets via nonlinear dynamical
systems. Second, we derive mathematical results concerning the behavior of 1D
PWL discontinuous maps that are of broader relevance, e.g., for applications
in economics, finance, engineering, biology, and ecology. The main
characteristic of our financial market model is that it corresponds to a 1D
discontinuous PWL map, defined over several partitions, with the real fixed
point being the same for the functions in all partitions. While such maps
produce dynamics that may be mistaken for chaotic behavior, our insights allow
us to conclude that they are, in fact, quasiperiodic.

To be precise, the financial market model we propose in the next section
belongs to a general class of maps, defined as follows: 

\begin{Definition} \textbf{(the class of 1D maps)}: We consider the class
of 1D PWL maps with a finite number of discontinuity points and functions
sharing the same real fixed point.
\end{Definition}

Before presenting our main results, formally stated in Theorem 1, a few
comments are in order. 

We always translate the unique real fixed point of the linear functions to the
origin, denoted as $O$, thereby making the system homogeneous. Note that if
the real fixed point belongs to a partition with slope equal to 1, so that it
is nonhyperbolic, we first consider the system with the hyperbolic fixed
point, translating it to $O$, and after it can be considered also the case of
nonhyperbolic fixed point. From now on, we assume that the fixed point is in
the origin. 

An essential element in the class of maps stated in Definition 1 is
discontinuity. The real fixed point $O$ may either be internal to one
partition or located at the boundary of two partitions. In the first case, the
map is continuous and differentiable at the fixed point $O$ (making it a
virtual fixed point for the functions in the other partitions). In the second
case, the map is continuous but piecewise smooth in the two partitions, with a
kink at the fixed point $O$. In this case, there must be at least one
additional discontinuity point to ensure a discontinuous map.

The map of the financial market model proposed in the next section satisfies
Definition 1, both when the real fixed point is internal to one partition and
when it is on the boundary of two partitions.

Despite being motivated by a financial market model, systems in the class of
maps defined above also appear in other application areas, such as electronic
engineering \cite{Kollar}, \cite{Zhu}, biology, and ecology \cite{Segura}.
Moreover, a map belonging to our class may appear as a restriction or first
return map on some straight-line segment of a two-dimensional PWL
discontinuous map, as proposed in \cite{GRSSW-24, GRSSW-24b,GRSSW-24c}. 

The one-dimensional bifurcation diagrams depicted in Fig.~1 illustrate the
dynamics of the map related to the financial market model proposed in the next
section. At first sight, the dynamics may appear chaotic. However, we will
show that this conclusion is not justified. Recall that for a discontinuous 1D
map $x^{\prime}=G(x)$, the most widely used definition of chaos is as
follows:

\begin{Definition} A map $G$ is chaotic in a closed
invariant set $X$ if periodic points are dense in $X$ and
there exists an aperiodic trajectory dense in $X$ (implying
topological transitivity).
\end{Definition} 

See, e.g., \cite{Devaney}. In this paper, "invariant" means mapped exactly
onto itself, $G(X)=X$.\footnote{In other papers, its meaning may be
also $G(X)\subseteq X,$ so that $X$ may be strictly mapped into itself.}

\begin{figure}
	\begin{center}
		\includegraphics[width=0.8\textwidth]
		{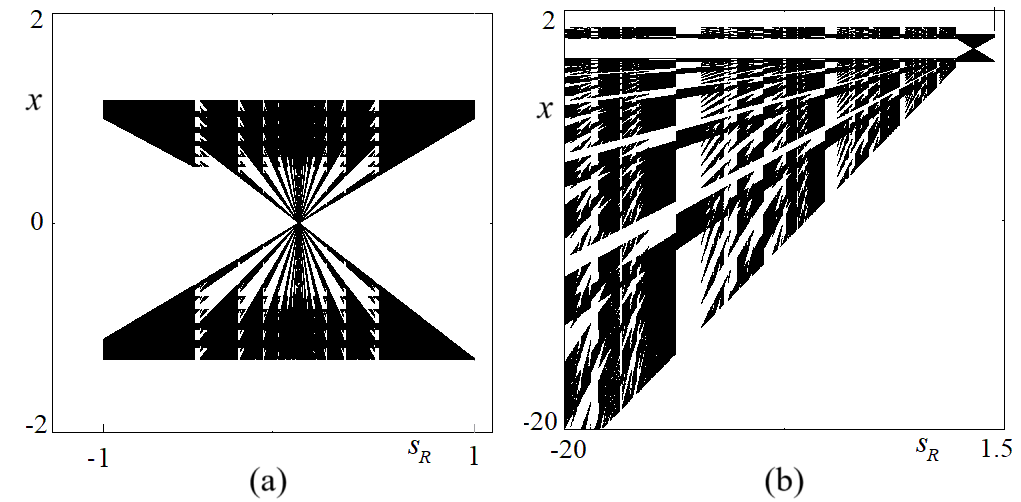} 
		\caption{\label{f1} One-dimensional bifurcation diagram of the map in (\ref{new})
			as a function of slope $s_{R}$. In (a), $Z_{L}=0.6,$ $Z_{R}=1,$ $s_{L}=-0.9,$
			$s_{M}=-1.3.$ In (b), $Z_{L}=1,$ $Z_{R}=0.6,$ $s_{L}=-0.5,$ $s_{M}=-1.1.$}
	\end{center}
\end{figure}

We show that chaos cannot emerge in this class of maps, since nonhyperbolic
cycles different from $O$ cannot exist. We prove that the dynamics of the maps
satisfying Definition 1 can be reduced to those of a PWL circle map
(\cite{deMelo}, \cite{Avrutin}). As a result, the only bounded dynamics
different from those related to the fixed point $O$ are either associated with
nonhyperbolic cycles or quasi-periodic trajectories densely filling certain
intervals. We prove the following theorem: 

\begin{Theorem} Let $x^{\prime}=G(x)$ be a 1D
	discontinuous PWL homogeneous map as in Definition 1. Then the following statements hold: 
	\begin{enumerate}
		\item A hyperbolic cycle different from the fixed point $O$ 
		cannot exist (and thus, there is no chaotic set).
		\item The only possible bounded invariant sets of map 
		$G,$ different from those related to the fixed point $O$ 
		(hyperbolic or nonhyperbolic), are those occurring in a PWL circle map. These
		consist of intervals dense with nonhyperbolic periodic points or quasiperiodic
		orbits dense in some intervals, with possible coexistence.
		\item Quasiperiodic orbits lead to (weak) sensitivity to initial
		conditions. 
		\item The Lyapunov exponent is zero.
	\end{enumerate}
\end{Theorem}

We continue as follows. In Sec.~2, we introduce the financial market model
that leads to a discontinuous map in three or four partitions. In Sec.~3, we
recall the properties of the PWL Lorenz map with one discontinuity to
emphasize that the PWL circle map separates regular regimes from chaotic ones.
The properties related to maps in our definition are considered in Sec.~4.
First, we prove the result in the simplest case with a unique discontinuity (Lemma 1). Second, we consider two discontinuity points, showing several
scenarios with the financial market model. Third, we consider the generic
case, proving Theorem 1. Fourth, we discuss the role of maps satisfying
Definition 1 as a separator of regular regimes from chaotic ones. Some
conclusions are drawn in Sec.~5.

\section{A simple nonlinear financial market model}

The seminal paper by Day and Huang [\cite{Day} has motivated numerous
follow-up studies aimed at explaining the intricate behavior of financial
markets (see \cite{Chiarella-09,Hommes-09,Lux-09,W-09} for surveys). In \cite{HD-93}, the authors developed a piecewise
linear version of their original financial market model. The simplicity of the
piecewise linear nature of their model provides new insights into the
properties of financial markets. Until recently, only a few piecewise linear
models had been proposed, mainly because the theory of PWL maps was not well
developed or widely known, at least in applied contexts. Fortunately, the
properties of 1D and 2D PWL maps have received greater attention in recent
years, leading to the emergence of many new models that enhance our
understanding of the dynamics of financial markets (see, e.g., \cite{H10,TWG-10,TGW-11,TWG-13}).

Following the approach in \cite{TWG-14}, we consider a financial market model
in which there are two types of chartists and two types of fundamentalists.
Type 1 chartists buy assets when they perceive a bull market and sell assets
when they perceive a bear market. Type 2 chartists follow the same trading
strategy; however, they become active only if the price deviates by at least a
certain distance from its fundamental value. Type 1 fundamentalists buy assets
when the market is undervalued and sell assets when it is overvalued. Type 2
fundamentalists do the same but---similar to type 2 chartists---enter the
market only when the price deviates by at least a certain distance from its
fundamental value.

Let $P_{t}$ be the log of the price in period $t$ and $F$ its log fundamental
value. The market maker mediates speculators' transactions out of equilibrium
and adjusts prices via excess demand using a log-linear price adjustment rule.
To be precise, the market maker quotes the log price for period $t+1$ as%
\begin{equation}
P_{t+1}=P_{t}+a(DC_{t}^{1}+DF_{t}^{1}+DC_{t}^{2}+DF_{t}^{2}) 
\label{MaMe}%
\end{equation}
The orders placed by type 1 chartists are formalized as%
\begin{equation}
DC_{t}^{1}=\left\{
\begin{array}
[c]{c}%
a^{1.b}+c^{1,b}(P_{t}-F)\\
-b^{1,d}+c^{1,d}(P_{t}-F)
\end{array}
\right.  \ \
\begin{array}
[c]{c}%
\text{for} \ \ P_{t}-F\geq0\\
\text{for} \ \ P_{t}-F<0
\end{array}
\label{C1}%
\end{equation}
Type 2 chartists are inactive when $-Z_{L}<P_{t}-F<Z_{R}$. Their orders are
expressed as%
\begin{equation}
DC_{t}^{2}=\left\{
\begin{array}
[c]{c}%
a^{2,b}+c^{2,b}(P_{t}-F)\\
0\\
-b^{2,d}+c^{2,d}(P_{t}-F)
\end{array}
\right.  \ \
\begin{array}
[c]{c}%
\text{for} \ \ P_{t}-F\geq Z_{R}\\
\text{for} \ \ -Z_{L}<P_{t}-F<Z_{R}\\
\text{for} \ \ P_{t}-F\leq-Z_{L}%
\end{array}
\label{C2}%
\end{equation}
The orders placed by type 1 fundamentalists are formalized as%
\begin{equation}
DF_{t}^{1}=\left\{
\begin{array}
[c]{c}%
-a^{1,b}+f^{1,b}(F-P_{t})\\
b^{1,d}+f^{1,d}(F-P_{t})
\end{array}
\right.  \ \
\begin{array}
[c]{c}%
\text{for} \ \ P_{t}-F\geq0\\
\text{for} \ \ P_{t}-F<0
\end{array}
\label{F1}%
\end{equation}
Type 2 fundamentalists are inactive when $-Z_{L}<P_{t}-F<Z_{R}$. Their orders
are given by
\begin{equation}
DF_{t}^{2}=\left\{
\begin{array}
[c]{c}%
-a^{2,b}+f^{2,b}(F-P_{t})\\
0\\
b^{2,d}+f^{2,d}(F-P_{t})
\end{array}
\right.  \ \
\begin{array}
[c]{c}%
\text{for} \ \ P_{t}-F\geq Z_{R}\\
\text{for} \ \ -Z_{L}<P_{t}-F<Z_{R}\\
\text{for} \ \ P_{t}-F\leq-Z_{L}%
\end{array}
\label{F2}%
\end{equation}
We assume that all reaction parameters in equations (\ref{C1}) to (\ref{F2})
are non-negative. Inserting (\ref{C1})-(\ref{F2}) into (\ref{MaMe}),
defining
\begin{equation}
	\begin{split}
s_{1}=a(c^{1,b}-f^{1,b})\ ,\ s_{2}=a(c^{1,d}-f^{1,d}),\\ s_{3}=a(c^{2,b}-f^{2,b})\ ,\ s_{4}=a(c^{2,d}-f^{2,d})
\end{split} 
\end{equation}
and expressing the log of the price in deviations from its log fundamental
value, i.e., $x_{t}=P_{t}-F$, we obtain the following system%
\begin{equation}
G_{4}:\left\{
\begin{array}
[c]{c}%
G_{L}:x^{\prime}=s_{L}x=(1+s_{2}+s_{4})x\ \ \ \ \\
G_{M-}:x^{\prime}=s_{M-}x=(1+s_{2})x\ \ \ \ \ \ \ \\
G_{M+}:x^{\prime}=s_{M+}x=(1+s_{1})x\ \ \ \ \ \ \ \\
G_{R}:x^{\prime}=s_{R}x=(1+s_{1}+s_{3})x\ \ \ \
\end{array}
\right.
\begin{array}
[c]{c}%
\text{for} \ x\leq-Z_{L}\\
\text{for} \ -Z_{L}<x\leq0\\
\text{for} \ 0\leq x<Z_{R}\\
\text{for} \ x\geq Z_{R}%
\end{array}
\label{gen}%
\end{equation}
The prime symbol stands for the unit time advancement operator. Moreover, the
parameters $s_{L}=(1+s_{2}+s_{4}),$ $s_{M-}=(1+s_{2}),$  $s_{M+}=(1+s_{1})$ and $s_{R}=(1+s_{1}+s_{3})$ represent the slopes of the map's four linear
functions. The four slope parameters may take any sign, positive or negative.
This is a map that satisfies Definition 1, with the fixed point located on the
boundary of two partitions and two additional discontinuity points.

The particular case where $Z_{L}=Z_{R},$ $s_{2}=s_{1}$ and $s_{4}=s_{3}$ has
been considered in \cite{TWG-15} and \cite{GT-2012}, where it has been shown
that the system can have cycles of any period that densely fill some segments
and quasiperiodic orbits that are dense in certain intervals, leading to
dynamics that may appear chaotic but are not.

Now let us relax some assumptions, assuming only that $s_{2}=s_{1}.$ We then
have the following system:%
\begin{equation}
G_{3}:\left\{
\begin{array}
[c]{c}%
G_{L}:x^{\prime}=s_{L}x=(1+s_{1}+s_{4})x\ \ \ \ \\
G_{M}:x^{\prime}=s_{M}x=(1+s_{1})x\ \ \ \ \ \ \ \ \ \ \\
G_{R}:x^{\prime}=s_{R}x=(1+s_{1}+s_{3})x\ \ \ \
\end{array}
\right.
\begin{array}
[c]{c}%
\text{for} \ x\leq-Z_{L}\\
\text{for} \ -Z_{L}<x<Z_{R}\\
\text{for} \ x\geq Z_{R}%
\end{array}
\label{new}%
\end{equation}
with $s_{M}=(1+s_{1}),$ $s_{L}=(s_{M}+s_{4})$ and $s_{R}=(s_{M}+s_{3}).$ This
map satisfies Definition 1 with the fixed point internal to the middle partition.

One of the threshold parameters $Z_{L}$ and $Z_{R}$ in (\ref{gen}) and
(\ref{new}) can be considered a scaling factor, and thus, be set at any value.
However, we keep them because they are economically relevant.

We also notice that when parameter $a^{1,b}$ in the trading rule of type 1
chartists differs to that in the trading rule of type 1 fundamentalists, then
the maps in (\ref{gen}) and (\ref{new}) have some offsets (i.e., the function
in some partition becomes affine instead of linear). The same applies for
parameters $a^{2,b}$, $b^{1,d}$ and $b^{2,d}$. In such cases, the systems do
not belong to the class of maps considered in this work. The PWL systems with
some affine functions can clearly have chaotic dynamics (see, e.g., the
systems in \cite{TGW-11,TWG-13,TWG-14}).

Differently, the dynamics that can occur in the system in (\ref{gen}), as well
as in the system in (\ref{new}), cannot be chaotic, as is the case with any
map satisfying Definition~1. However, from the applied point of view, the
systems may nevertheless mimic important properties of financial market dynamics.

\section{Preliminary properties of a PWL Lorenz map}

In this section, we recall some properties of the dynamics of 1D PWL
discontinuous maps with a single discontinuity point, known as 1D PWL Lorenz
maps, and circle maps.

Any map $f:[a,b]\rightarrow\lbrack a,b]$ such that $f(a)=f(b)\in(a,b)$ with a
unique discontinuity point $c\in(a,b)$ where $f=f_{L}$ is monotonically
increasing in $[a,c)$ with $f_{L}(c)=b$ and $f=f_{R}$ is monotonically
increasing in $(c,b]$ with $f_{R}(c)=a$, defines a homeomorphism of the circle
\cite{deMelo}. Note that then $f_{L}(f_{R}(c))=f_{R}(f_{L}(c)).$

A \textbf{PWL circle map} is a circle homeomorphism in which the two
components are linear or affine functions.

A 1D \textbf{PWL Lorenz map} may be defined as follows:
\begin{equation}
f=\left\{
\begin{array}
[c]{c}%
f_{L}:x^{\prime}=a_{L}x+\mu_{L}\ \ \ \text {for} \ x<0\\
f_{R}:x^{\prime}=a_{R}x+\mu_{R}\ \ \ \text{for} \ x>0
\end{array}
\right.  \label{Lorenz}%
\end{equation} 
with $a_{R/L}>0$ and $\mu_{L}\neq\mu_{R}$ (so that $f$ is discontinuous)$.$
For this class of maps, assuming $\mu_{R}<\mu_{L},$ and the interval
$I=[f_{R}(0),f_{L}(0)]=[\mu_{R},\mu_{L}]$ an absorbing interval (that is, a
point in a neighbourhood of $I$ enters $I$\ in a finite number of iterations,
and $f(I)\subseteq I$).

\begin{figure}[ht]
	\begin{center}	
		\includegraphics[width=.8\textwidth]
		{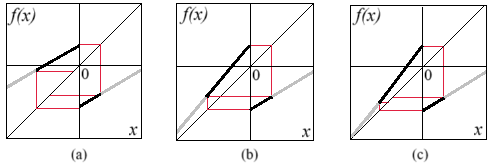} 
		\caption{\label{f2} PWL Lorenz map. In (a), a gap map; in (b), a circle map; in (c), an
			overlapping map.}
	\end{center}	
\end{figure}

Considering the images of the discontinuity point $x=0,$ the dynamics of a PWL
Lorenz map, as depicted in Fig.~2, can be classified into three different types \cite{Berry,Rand}: 

(L1) In the case $f_{R}\circ f_{L}(0)<f_{L}\circ f_{R}(0),$ the map is called
a \textbf{gap map}, where the gap is the interval $J=]f_{R}\circ
f_{L}(0),f_{L}\circ f_{R}(0)[$ and $f(I)\subset I.$ Each point $x\in I$ has
either a unique preimage or no preimages inside $I$. The rotation number is
well defined (see \cite{deMelo}); it is the same for any point in $I$, and its
dynamics can only be as follows:

- If the rotation number of map $f$ is rational, then the $\omega-$limit set
for all $x\in I$ is an attracting hyperbolic cycle.

- If the rotation number is irrational, then the $\omega-$limit set for all
$x\in I$ is a unique Cantor set attractor $\mathcal{C}$,
\[
\mathcal{C}=\underset{n\geq0}{\cap}f^{n}(I)=I\backslash\underset{n\geq0}{\cup
}f^{n}(J)
\]
The Cantor set attractor $\mathcal{C}$ is the closure of quasiperiodic orbits.
A Cantor set attractor has a self-similar structure, as is the case with any
Cantor set.

In this class of maps, the irrational rotation is the nongeneric case, which
is not persistent under parameter perturbation, while the rational rotations
occur in open regions of the parameter plane and are thus persistent under
parameter perturbations. The existence regions related to cycles are bounded
by border collision bifurcation curves following a period adding rule and
Farey summation rule \cite{Avrutin}.  

(L2) Map $f$ is a \textbf{circle map} if and only if the following equality
holds: $f_{R}\circ f_{L}(0)=f_{L}\circ f_{R}(0)$. Now $I$ is invariant,
$f(I)=I.$ In this case, we have a circle homeomorphism. It is known
\cite{deMelo} that the rotation number is well defined; it is the same for any
point of interval $I$, and the dynamics depend on the rotation number as follows:

- If the rotation number of map $f$ is rational, then the absorbing interval
$I=[f_{R}(0),f_{L}(0)]$ is filled with periodic points (of the same period).

- If the rotation number of map $f$ is irrational, then the absorbing interval
$I$ is filled with quasiperiodic orbits dense in that interval (i.e., any
orbit of $f$ is quasiperiodic and dense in $I$, so that the $\omega$-limit set
of each orbit is the whole interval).

The irrational rotation is now the generic case, while the case of a rational
rotation is not persistent under parameter perturbations (that is, in the
parameter space satisfying $f_{R}\circ f_{L}(0)=f_{L}\circ f_{R}(0),$ the
rotation number is irrational with a probability of 1, since the rational
rotation occurs in a parameter set of zero Lebesgue measure).

(L3) In the case $f_{R}\circ f_{L}(0)>f_{L}\circ f_{R}(0),$ the map is called
an \textbf{overlapping map}, since the interval $[f_{L}\circ f_{R}%
(0),f_{R}\circ f_{L}(0)]$ consists of points having two different rank-one
preimages in the invariant interval $I=[f_{R}(0),f_{L}(0)],$ $f(I)=I.$ This
system can have different kinds of dynamics, especially hyperbolic cycles and
chaotic attracting or repelling sets (which necessarily imply the existence of
infinitely many unstable cycles). In the PWL case of interest here, this
regime is characterized by attracting chaotic intervals, bounded by images of
the discontinuity point.

Thus, the case of a circle map denotes the transition from regular dynamics
(without chaos, in a gap map) to chaotic dynamics (existence of chaotic sets,
in the overlapping map), see \cite{G-2010,Avrutin}. 

\section{1D discontinuous PWL homogeneous maps}

In this section, we study the properties of maps satisfying Definition 1. In
Section 4.1, we consider maps with a unique discontinuity point. In Section
4.2, we examine maps with two discontinuity points. In Section 4.3, we prove
Theorem 1. In Section 4.4, we discuss the role of maps satisfying Definition 1
as separators between regular regimes and chaotic ones. The last two
subsections also include comments on the implications of our financial market model.

\subsection{One discontinuity point}

In the simplest case, the map consists of only two partitions and is defined
as%
\begin{equation}
F=\left\{
\begin{array}
[c]{c}%
F_{L}:x^{\prime}=s_{L}x\ \ \ \text{for} \ x<h\\
F_{R}:x^{\prime}=s_{R}x\ \ \ \text{for} \ x>h
\end{array}
\ \ h\neq0\right.  
\label{MapF}
\end{equation}
In order to have a unique fixed point in the origin and a unique discontinuity
point, it must be $h\neq0$.

In the following lemma, we prove that besides the possible occurrence of
divergent trajectories and trajectories related to the fixed point $O,$ the
only possible bounded dynamics are those occurring in a PWL circle
map. 

\begin{Lemma} Let $F$ be the map defined in (\ref{MapF}).
Then the only possible bounded $\omega$-limit sets different from the
fixed point $O,$ and from nonhyperbolic cycles associated with the
nonhyperbolic fixed point $O$, are those of a PWL circle map.
\end{Lemma}

Proof. To fix the reasoning, let us assume $h>0$ (which can also serve as a
scaling factor, and be set to $1$). The case $h<0$ is topologically conjugate
to $h>0$.

When one of the two slopes is zero, then the only possible bounded 
$\omega$-limit set is the fixed point $O$. So, let us consider $s_{L,R}\neq0$.
Moreover, if $|s_{L}|=1$ and/or $|s_{R}|=1,$ the bounded $\omega-$limit sets
different from the fixed point $O$ can only be nonhyperbolic fixed points or
nonhyperbolic 2-cycles (filling densely some segments).

Let us now consider the generic case, proving the statement for the different
signs of the slopes of the linear functions.

(i1) In the case of \textbf{positive slopes}, asymptotic dynamics different
from the fixed point $O$ or divergence can occur only for $F_{L}%
(h)>h>F_{R}(h),$\ leading to $s_{L}>1>s_{R},$ and in this case it is
$F_{R}\circ F_{L}(h)=s_{R}s_{L}h=F_{L}\circ F_{R}(h).$ Thus, we obtain a PWL
circle map in the absorbing interval $I=[F_{R}(h),F_{L}(h)]=[s_{R}h,s_{L}h].$

(i2) In the case of \textbf{one negative and one positive slope}, asymptotic
dynamics different from the fixed point $O$ or divergence can occur only for
$s_{L}<0$ and $s_{R}>0,$ with $F_{L}\circ F_{L}(h)>h>F_{R}(h),$\ leading to
$s_{L}^{2}>1>s_{R}.$ In this case, the dynamics can be obtained via the first
return map in $x>0$, that is given by $F_{L}\circ F_{L}(h)$ for $x<h$ and
$F_{R}(h)$ for $x>h,$ so that 
$F_{R}\circ F_{L}^{2}(h)=s_{R}s_{L}^{2}h=F_{L}^{2}\circ F_{R}(h).$ 
Thus, we obtain a PWL circle map in the absorbing
interval $I=[F_{R}(h),F_{L}^{2}(h)]=[s_{R}h,s_{L}^{2}h].$

(i3) In the case of negative slopes, the dynamics can be studied via the
second iterate, and asymptotic dynamics different from the fixed point $O$ or
divergence can occur only for 
$F_{L}\circ F_{L}(h)>h>F_{L}\circ F_{R}(h),$ leading to $s_{L}^{2}>1>s_{L}s_{R}.$ 
In this case, its dynamics are the
same as those of the first return map, in partition $x>0$, or in $x<0,$ in a
neighborhood of the discontinuity point. In the partition $x>0,$ it is given
by $F_{L}\circ F_{L}(h)$ for $x<h$ and $F_{L}\circ F_{R}(h)$ for $x>h,$ so
that $F_{L}\circ F_{R}\circ F_{L}^{2}(h)=s_{R}s_{L}^{3}h=F_{L}^{3}\circ
F_{R}(h),$ and it is a PWL circle map in the absorbing interval 
$I=[F_{L}\circ F_{R}(h),F_{L}\circ F_{L}(h)]=[s_{L}s_{R}h,s_{L}^{2}h].$ In the partition
$x<0$, it is given by $F_{R}\circ F_{L}(h)$ for $x<F_{L}^{-1}(h)=h/s_{L}$ and
$F_{L}\circ F_{L}(h)$ for $x>F_{L}^{-1}(h)=h/s_{L},$ so that 
$F_{L}\circ F_{L}\circ F_{R}\circ F_{L}(h/s_{L})=s_{L}^{2}s_{R}h=F_{R}\circ F_{L}^{3}(h/s_{L}),$ 
and it is a PWL circle map in the absorbing interval
$J=[F_{L}\circ F_{L}(h/s_{L}),F_{R}\circ F_{L}(h/s_{L})]=[s_{L}h,z_{L}s_{R}h].$ 
Both lead to the same dynamics of map $F.\square$ 
\smallskip
\begin{figure}[ht]
	\begin{center}	
		\includegraphics[width=1\textwidth]{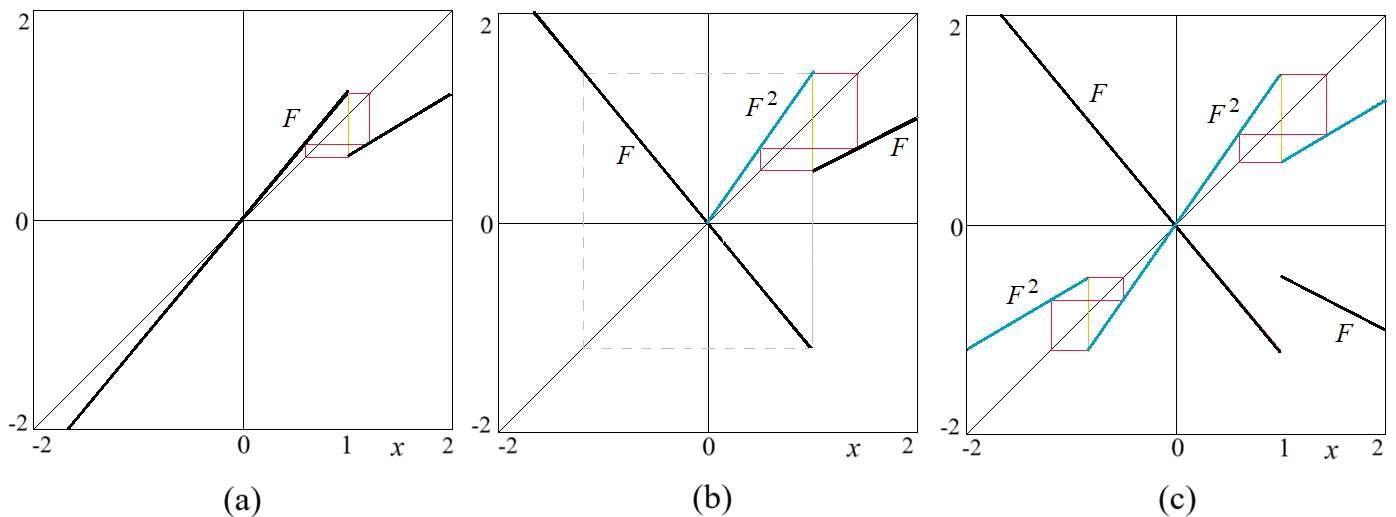} 
		\caption{\label{f3} Map $F$ in (\ref{MapF}) with $h=1$. Function $F$ is shown in black and
			its second iterate in azure. In (a), case (i1); in (b), case (i2); in (c),
			case (i3).}
	\end{center}
\end{figure}

The three possible generic cases described in the proof of Lemma 1 are shown
in Fig.~3.

\textbf{Remarks}. It is worth noting that, unlike the dynamics occurring in a
gap map, and similarly to the comment in the case (L2) of a Lorenz map in the
previous section (PWL circle map), the generic case is an irrational rotation
number, and thus, the trajectories are quasiperiodic and dense in intervals.

For map $F$, the existence of dense cycles is associated with a set of zero
Lebesgue measure in the $(s_{R},s_{L})$ parameter plane of the possible
slopes$.$ In the case of positive slopes, cycles occur if and only if integers
$p$ and $q$ exist such that $s_{L}^{p}s_{R}^{q}=1.$ This leads to
nonhyperbolic cycles (filling densely interval $I$) of period $n=p+q,$
assuming $p$ and $q$ are the smallest integers, and rotation number $\rho=p/n$
(or $(1-\rho)=q/n$).

It is also worth noting that the map known as "rigid rotation", that is, $R:$
$x\rightarrow x+\rho\ (\operatorname{mod}1),$ for $\rho\in(0,1),$ is not
included in our definition, since considering the map $x^{\prime}=x+\rho$ for
$0\leq x<(1-\rho),$ and $x^{\prime}=x+\rho-1$ for $(1-\rho)<x\leq1,$ the fixed
point is at infinity. However, this circle homeomorphism with rotation number
$\rho$ is topologically conjugate to map $F$ with positive slopes and the same
rotation number.

\subsection{Two discontinuity points}

Let us comment on the case of two discontinuity points. 
For a 1D map satisfying Definition 1 with two discontinuities and the fixed
point internal to the middle partition, the map is given, for $h_{1}<h_{2}$,
$h_{1}\neq0,$ $h_{2}\neq0$, by:%
\begin{equation}
G_{1}:\left\{
\begin{array}
[c]{c}%
G_{L}:x^{\prime}=s_{L}x\ \ \ \ \text{for} \ x<h_{1}\ \ \ \ \ \ \ \\
G_{M}:x^{\prime}=s_{M}x\ \ \ \text{for} \ h_{1}<x<h_{2}\\
G_{R}:x^{\prime}=s_{R}x\ \ \ \text{for} \ x>h_{2}\ \ \ \ \ \ \
\end{array}
\right.  
\label{mapG}%
\end{equation}
A particular case of map (\ref{mapG}) with $h_{1}=-h_{2}$ and $s_{R}=s_{L}$
has been considered in \cite{GT-2012} and \cite{TWG-15}. In these works, it
has been shown that the only bounded trajectories, different from the fixed
point $O,$ and the cycles related to the fixed point $O$ when it is
nonhyperbolic, are those occurring in a PWL circle map.

In the generic case, the same result can be obtained. A few examples with map
$G_{1}$ in (\ref{mapG}), which can be considered as our applied model in
(\ref{new}), with a different denomination of the discontinuity points, are
shown in Figs.~4 and 5, where it is evidenced that it is sufficient to look
for the first return map in a suitable interval (sometimes in a neighborhood
of a discontinuity point, but not necessarily), to see that the first return
map is a PWL circle map. The fact that coexistence of attractors is possible
can be seen also from the function in azure in Fig.~3(c) (although in that
figure it is the second iterate of map $F$). 

\begin{figure}[ht]
	\begin{center}	
		\includegraphics[width=.7\textwidth]{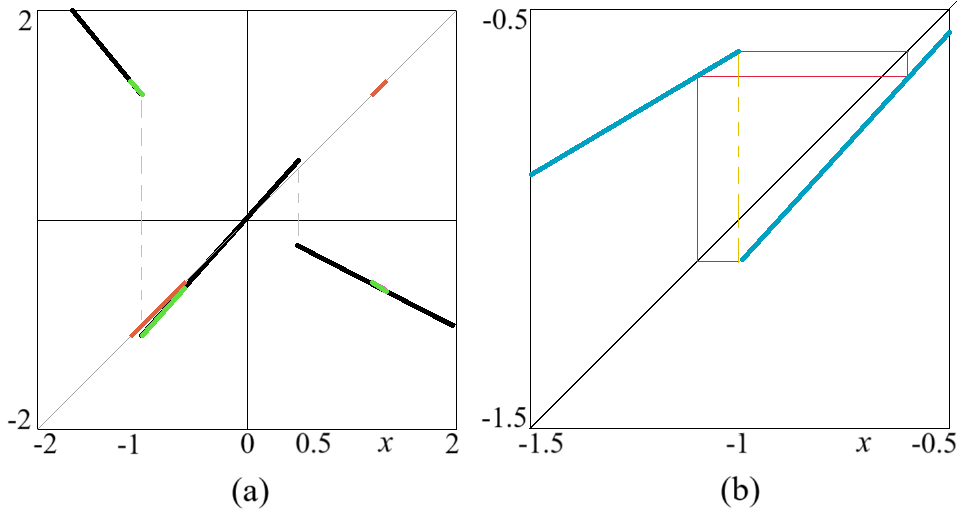} 
		\caption{\label{f4} Map $G_{3}$ in (\ref{new}) with $Z_{L}=1,$ $Z_{R}=0.5$, $s_{L}=-1.2,$
			$s_{M}=1.2,$ $s_{R}=-1.1.$ In (b), first return map in a neighborhood of
			$Z_{L}=1$.}
	\end{center}
\end{figure}
\begin{figure}[ht]
	\begin{center}	
		\includegraphics[width=0.7\textwidth]
		{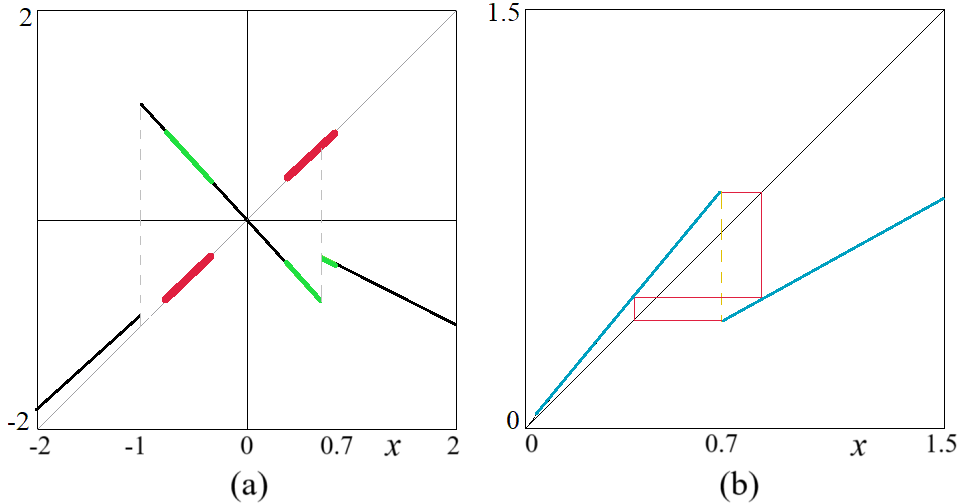} 
		\caption{\label{f5} Map $G_{3}$ in (\ref{new}) with $Z_{L}=1,$ $Z_{R}=0.7,$ $s_{L}=0.9,$
			$s_{M}=-1.1,$ $s_{R}=-0.5.$\ In (b), first return map in a neighborhood of
			$Z_{R}=0.7$.}
	\end{center}
\end{figure} 

To get the first return map in some interval, each point $x$ is iterated and
when the trajectory returns for the first time to the interval, we assign to
$x$ that value, say $F_{ra}(x).$ Considering the symbolic sequences of the
trajectory in the usual way, using the partitions involved in the trajectory,
the first return has a fixed symbolic sequence$.$ Due to the linearity, the
definition will be the same for an interval of points (i.e., the same symbolic
sequence), up to a point, say $x=c,$ whose trajectory merges with a
discontinuity point of the map. For $x>c,$ the first return leads to a
different function $F_{rb}(x)$, i.e., a different symbolic sequence. Since the
composition of linear homogeneous functions is always a linear homogeneous
function, the branches are functions with the same fixed point in the origin.
Then we have a map such as map $F$ in (\ref{MapF}), and the dynamics can be
ultimately described using $x^{\prime}=F_{ra}(x)$ for $x<c$ and 
$x^{\prime}=F_{rb}(x)$ for $x>c$ in the considered interval.

In the examples in Figs.~4 and 5, we have shown the first return map in an
interval including a discontinuity point. The map shown in Fig.~6(a) has many
intervals filled with quasiperiodic trajectories. The first return map in one
segment is shown in Fig.~6(b). Considering its discontinuity point $c$, the
rank-1 preimage $c_{-1},$ and the image $c_{1},$ we have two segments on the
right side of the discontinuity point $c$. Each point in the segment
($c,c_{-1}$) has the trajectory $R(LM)^{3}RM^{2}(ML)^{3}$, and then it returns
to the right side of $c$, while each point in segment ($c_{-1},c_{1}$) has the
trajectory $R(LM)^{3}R(LM)^{3}RM^{2}(ML)^{3}$, and then it returns to the
right side of $c.$

\begin{figure}[ht]
	\begin{center}	
		\includegraphics[width=0.9\textwidth]{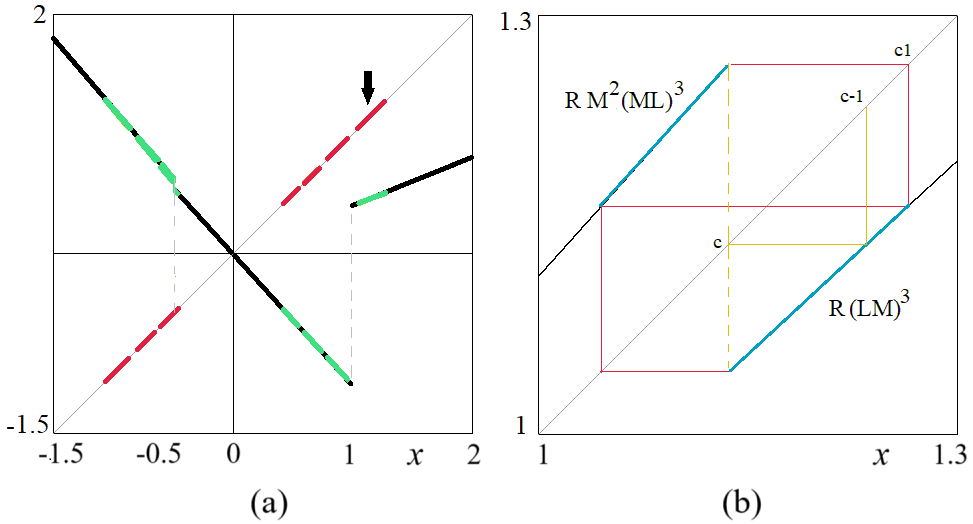} 
		\caption{\label{f6} Map $G_{3}$ in (\ref{new}) with $Z_{L}=0.5,$ $Z_{R}=1,$ $s_{L}=-1.2,$
			$s_{M}=-1.1,$ $s_{R}=-0.4.$ In (b), first return map in the rightmost
			interval, evidenced by an arrow in (a).}
	\end{center}
\end{figure}

In the case that the product of the slopes of the symbolic sequence
$R(LM)^{3}$ $RM^{2}(ML)^{3}$ is $1$, then the product of the two slopes in the
first return map is also $1$ (so that it is filled with 2-cycles) and map
$G_{3}$ has segments filled with periodic points of cycles of period 16.

In the example shown in Fig.~7, we consider the first return in a small segment
evidenced by an arrow in Fig.~7(a), which is shown in Fig.~7(b). This first
return map takes 5 iterations for $x<c$ with a symbolic sequence, say 
$S_{ra}$, and 30 iterations for $x>c,$ with a symbolic sequence, say $S_{rb}.$
Looking at the attractor in Fig.~7(a), it seems chaotic in several intervals,
but it is not, and its dynamics are predictable. Taking the preimages of the
discontinuity point $c$ of the first return map, the segment on the right side
is separated into two parts. Each point in the first segment, close to $c,$
has the trajectory $S_{rb}(S_{ra})^{5}$ and returns in the right side of $c$,
while each point in the second segment at the right side of $c,$ has the
trajectory $S_{rb}(S_{ra})^{4},$ and is back again in the right side of $c$.
Clearly, $S_{rb}$ and $S_{ra}$ are to be expanded with the 30 and 5 symbols,
respectively (related to the partitions $L$, $M$, $R$), $S_{rb}=R^{2}%
ML^{2}(MRML^{5})^{3}M,$ and $S_{ra}=R^{2}MLM$. Figure 7(c) shows a 1D
bifurcation diagram as a function of slope $s_{M}$ in the middle partition.

\begin{figure}[ht]
	\begin{center}	
		\includegraphics[width=1\textwidth]{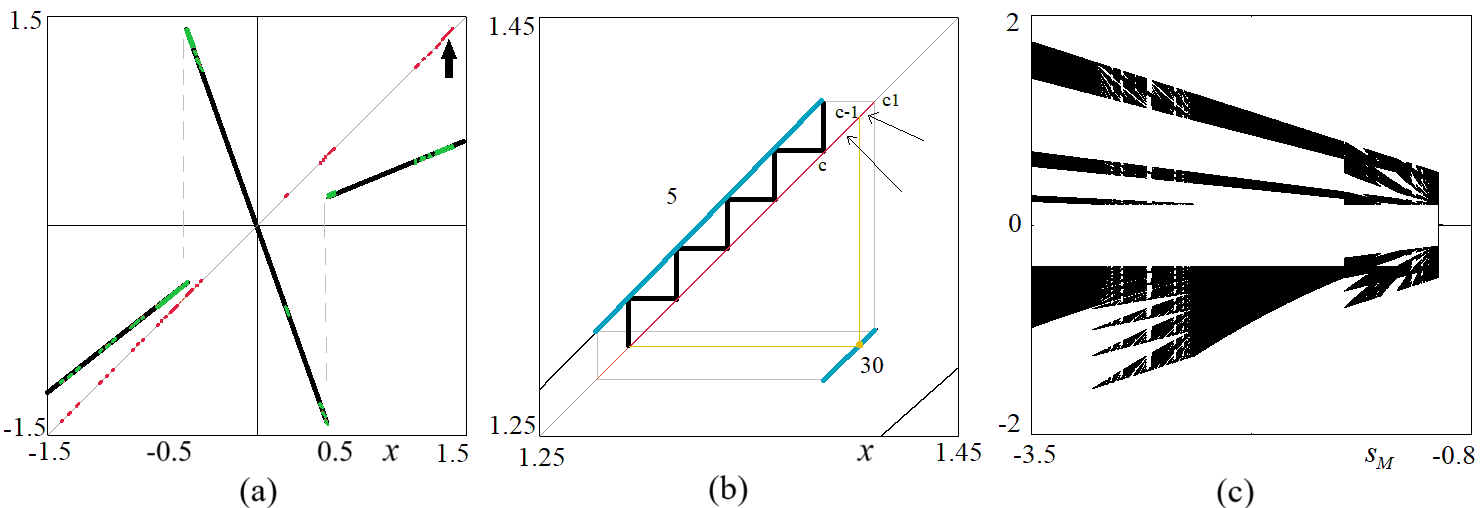} 
		\caption{\label{f7} Map $G_{3}$ in (\ref{new}) with $Z_{L}=0.5,$ $Z_{R}=0.5,$ $s_{L}=0.8,$
			$s_{M}=-2.82,$ $s_{R}=0.4.$ In (b), first return map in the rightmost
			interval, evidenced by an arrow in (a). In (c), 1D bifurcation diagram as a
			function of slope $s_{M}.$ }
	\end{center}
\end{figure}

A similar behavior always occurs, regardless of the number of partitions in
the map. In this kind of maps, hyperbolic cycles different from the fixed
point $O$ cannot exist, so that homoclinic cycles cannot exist either (and
thus there are no chaotic sets). The absence of hyperbolic cycles leads to the
dynamics of a map of Lorenz type, but it cannot be a gap map because a gap map
can have only hyperbolic cycles or Cantor set attractors, which are not
possible here. Neither it is possible to get an overlapping map (which in the
linear case has hyperbolic cycles different from the fixed point $O$, while
they cannot exist here). It follows that the only possible bounded dynamics
not converging to the fixed point $O$ are segments filled with nonhyperbolic
cycles (all with the same period) or quasiperiodic trajectories dense in some
intervals, as occurring in a PWL circle map.

\subsection{Proof of Theorem 1}
We are now prepared to prove Theorem 1, which characterizes the attracting
sets that may arise in 1D PWL maps, as defined in Definition 1, assuming, for
simplicity (without loss of generality), that the unique fixed point of the
functions defined in the different partitions has been set in the
origin. \medskip 

\textbf{Proof of Theorem 1.} Let us first consider maps with the fixed point
internal to a partition.

(1) Assuming, by contradiction, that a map in our class can have a hyperbolic
cycle of some period, with a given symbolic sequence, the periodic points must
be fixed points of the related composite map (obtained following the symbolic
sequence of the cycle). However, any composition of the linear maps in the
definition is also a linear map that has only the fixed point $O$ as its
possible hyperbolic fixed point. Thus, a hyperbolic cycle different from the
fixed point $O$ cannot exist, and obviously nor can a chaotic set, since one
of the properties of a chaotic set is that homoclinic hyperbolic cycles must
be dense in it.

(2) Since the fixed point is actual in one partition, it may be hyperbolic or
not. If not (the eigenvalue is $+1$ or $-1$), then a segment exists filled
with nonhyperbolic cycles (fixed points or 2-cycles), all with eigenvalue $+1$.

If bounded $\omega$-limit sets different from those related to the fixed point
$O$ (hyperbolic or nonhyperbolic) exist, these must involve different
partitions of the map (at least two), that is, at least one discontinuity
point. Let $A$ ($\neq O$) be a bounded $\omega$-limit set of a trajectory of
map $G$, then $G(A)=A$. If it is a $k$-cycle, it cannot be hyperbolic. Thus,
it is nonhyperbolic, and nonhyperbolic cycles in PWL maps in our definition
must be dense in some intervals. In fact, any iterate with map $G^{k}$ of an
interval $I$ in a partition is such that $G^{k}(I)$ consists of segments
belonging to straight lines issuing from the fixed point $O$, and thus cannot
have a kink point with a contact on the main diagonal (nor can transversely
intersect it). It follows that it can only be the case that the $k-$th iterate
of the map, $G^{k}$, includes segments on the main diagonal.

Note that any composition of the linear functions defining map $G$ consists of
linear branches, so that the first return map in any segment inside the
absorbing interval is always a function with the same characteristics as map
$G,$ the branches have the same fixed point $O$. Hence the first return map in
a suitable segment leading to a function with only one discontinuity point is
necessarily as map $F$ in (\ref{MapF}), and by Lemma 1 it can be only a PWL
circle map.

(3) Recall that a map $F:I\rightarrow I$ is said to have sensitive
dependence on initial conditions on $I$ if there exists a constant $\delta>0$
such that, for any $x\in I$ and any neighborhood $U(x)$ of $x$, there exists
an $x\in U(x)$ and a positive integer $n$ such that 
$|F^{n}(x)-F^{n}(x)|>\delta$ (see, e.g., \cite{Devaney,Wiggins}). In our case, when
there are quasiperiodic trajectories, these can be ultimately reduced to those
of a PWL circle map, since the first return map in some suitable interval
leads to a map as $F$ in (\ref{MapF}). Considering any two nearby points in
the first return PWL circle map, in a finite number of iterations, the
trajectories will have points separated by the discontinuity point. Thus, the
definition cannot be satisfied. However, after one iteration, the trajectories
will again be close to each other, leading to the property of weak sensitivity
(a proof for a map as map $F$ in (\ref{MapF}) is also presented in
\cite{Kollar}).

(4) Since the dynamics can ultimately be reduced to those of a PWL circle map
as map $F$ in (\ref{MapF}), it is well known that it has Lyapunov exponent
zero (a proof is also given in \cite{Kollar}).

For maps in which the fixed point $O$ is a kink point, the proof is the same
except for the comments related to the fixed point at the beginning of point
(2). In fact, since the map is now continuous at the origin, but not
differentiable, the behavior of the points close to the fixed point $O$ on
its right and left sides may be different, depending on the related slopes.
$\square$ 
\medskip

Some comments are in order. The example with map $G_{3}$ in (\ref{new}) shown
in Fig.~8 evidences that the first return map $F$ has an interval filled with
nonhyperbolic 2-cycles. Such cycles for map $G_{3}$ correspond to
nonhyperbolic 3-cycles with symbolic sequence $RLM,$ since it is 
$s_{L}s_{M}s_{R}=1$. Figure 8(c) shows, besides $G_{3},$ the third iterate
$(G_{3})^{3},$ which has segments on the main diagonal filled with
nonhyperbolic fixed points.

\begin{figure}[ht]
	\begin{center}	
		\includegraphics[width=1\textwidth]{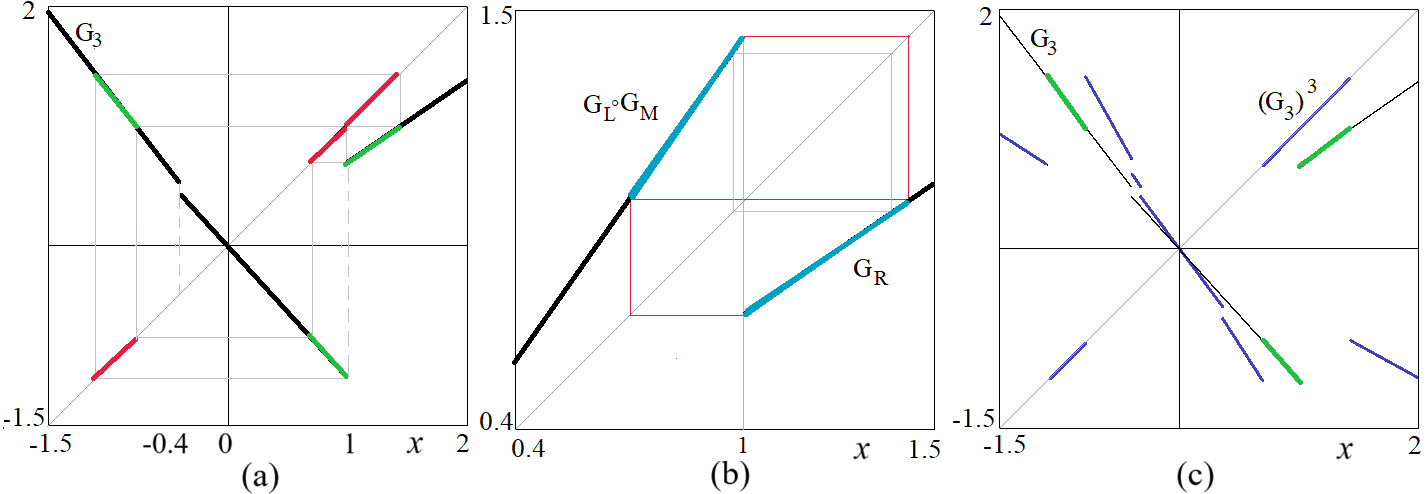} 
		\caption{\label{f8} Map $G_{3}$ in (\ref{new}) with $Z_{L}=0.4,$ $Z_{R}=1,$ $s_{L}=-1.3,$
			$s_{M}=-1.1,$ $s_{R}=0.6993006993$. In this example, the map has segments
			filled with cycles of period $3$ since $s_{L}s_{M}s_{R}=1.$ In (b), first
			return map in a neighborhood of $Z_{R}=1$. In (c), the map and its third
			iterate, having segments on the diagonal filled with fixed points.}
	\end{center}
\end{figure}
\begin{figure}[ht]
	\begin{center}	
		\includegraphics[width=1\textwidth]{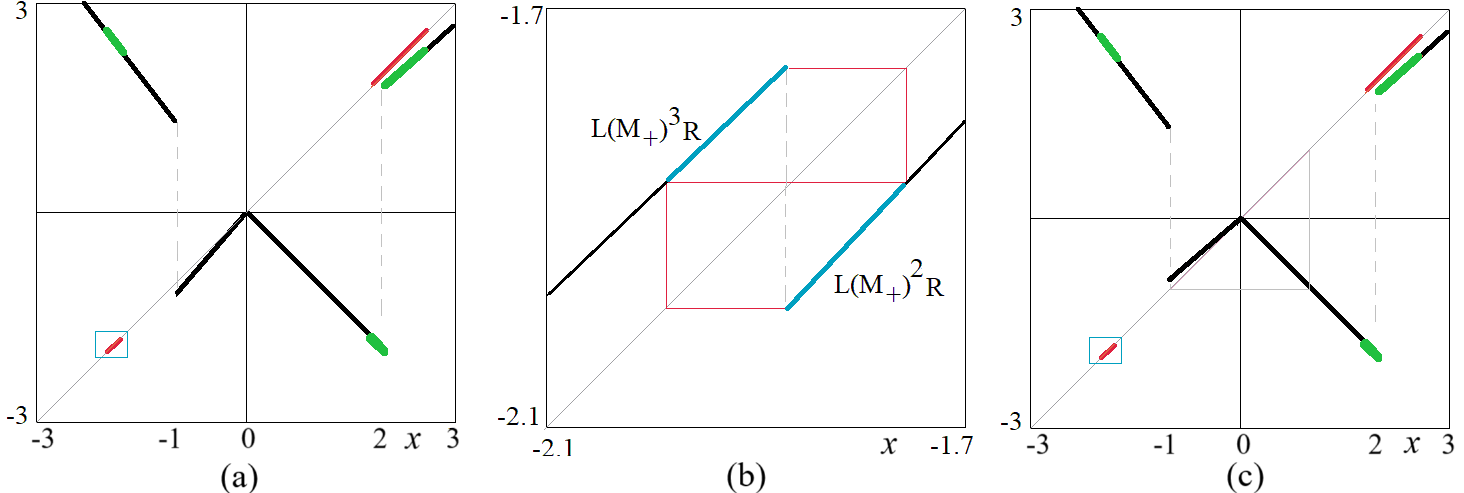} 
		\caption{\label{f9} Map $G_{4}$ (\ref{gen}) with$\ Z_{L}=1,$ $Z_{R}=2,$ $s_{L}=-1.3,$
			$s_{M-}=1.2,$ $s_{M+}=-1$, $s_{R}=0.9.$ In (b), first return map in the
			leftmost interval, evidenced by a rectangle in (a). In (c), the slopes of the
			function in the\ middle left partition is $s_{M-}=0.9.$ }
	\end{center}
\end{figure}

The example with map $G_{4}$ in (\ref{gen}) shown in Fig.~9 evidences that the
two partitions with the kink point in the fixed point $O$ may not be involved
in the asymptotic dynamics. Here, a right neighborhood of the origin is mapped
to the left side, and all the points in the left partition of the fixed point
$O$ are converging to an attractor belonging to the other partitions, as shown
via the symbolic sequence occurring in the first return map (in the interval
in the square of Fig.~9(a)). In this example, the slopes of the function in the
left/right partitions of the fixed point $O$ are $s_{M-}=1.2$ and $s_{M+}=-1,$
respectively.

However, changing the slope in the left partition $M_{-}$ of the fixed point
$O$, as in the example in Fig.~9(c), we have coexistence: the same intervals
filled with quasiperiodic trajectories shown in Fig.~9(a) coexist with the
attracting fixed point $O$, whose basin of attraction is the segment $(-1,1)$.
The remaining points converge to the other attractor. Now the slopes of the
function in the kink point $O$ are $s_{M-}=0.9$ and $s_{M+}=-1.$ 

It is worth noting that the invariant set $A$ that can occur in the class of
1D PWL maps considered here, different from the fixed point $O,$ can be an
attracting set or not (depending on the behavior of the map in points outside
the invariant set). The fixed point $O$ or the invariant segment related to
the nonhyperbolic fixed point $O$ may be attractors in Milnor
sense\ (attracting a set of positive Lebesgue measure, \cite{Milnor}) and may
coexist with intervals filled with quasiperiodic trajectories or periodic
points. Note that when the map has segments filled with nonhyperbolic cycles,
the segments lead to an attracting set, but not to an attractor, since such an
invariant set does not have a dense orbit. It follows that the only attractor
$A$ (an invariant set with a dense orbit and a suitable neighborhood $U(A)$
such that $A=\cap_{n>0}G^{n}(U(A))$) different from the fixed point $O$ may be
the set of intervals filled with quasiperiodic trajectories.

Moreover, we remark that in the parameter space of the slopes of a system in
our definition, the occurrence of segments filled with nonhyperbolic cycles is
not structurally stable. This can occur only when the product of the related
slopes (of the linear functions in the partitions with periodic points) leads
exactly to the value $1$ (as in the example shown in Fig.~8). Clearly, the
generic case is the occurrence of quasiperiodic trajectories. This is also the
reason why the numerical computations, as one parameter is varied, lead to
figures that seem related to chaos (as in Fig.~1(a,b) and Fig.~7(c)).

The quasiperiodic character of a trajectory may be better observed when an
applied system is investigated by using the versus-time trajectory, as the
prices in the financial market model. A couple of examples with the case shown
in Fig.~7 at $s_{M}=-1.5$ and $s_{M}=-2.82$ are reported in Fig.~10. 

\begin{figure}[ht]
	\begin{center}	
		\includegraphics[width=0.8\textwidth]{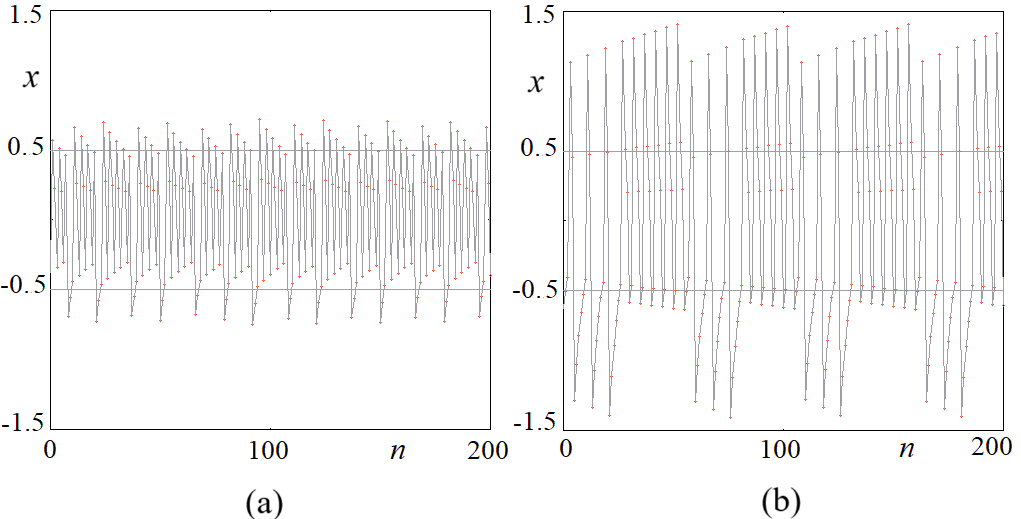} 
		\caption{\label{f10} Map $G_{3}$ in (\ref{new}) with $Z_{L}=0.5,$ $Z_{R}=0.5,$ $s_{L}=0.8,$
			$s_{R}=0.4.$ Versus time trajectory, for 200 iterations. In (a), $s_{M}=-1.5.$
			In (b), $s_{M}=-2.82.$}
	\end{center}
\end{figure}

Recall that $x$ represents the log price of the financial market in deviation
from its log fundamental value. Hence, when $x>0$, the financial market is
overvalued, and when $x<0$, the financial market is undervalued. Clearly, our
financial market model is capable of producing endogenous price dynamics.
Depending on the parameter setting, prices fluctuate around their fundamental
value within a narrower (left) or broader band (right), qualitatively
mimicking the excessively volatile boom-bust behavior of actual financial
markets. The dynamics depicted in the right panel, Fig.~10(b), are particularly
interesting in this respect. The financial market model more frequently
produces a bull market than a bear market, though when a bear market emerges,
it may be quite significant.

\subsection{Class of maps separating regular/chaotic behaviors}

The properties of the maps in our definition, described in Theorem 1, are due
to the existence of a unique fixed point for all the functions in the
partition. We may expect that when this assumption is not satisfied, then
hyperbolic cycles different from the fixed point $O$ may appear as well as chaos.

As for the PWL circle map in the Lorenz family recalled in Sec.~3, which
separates regular regimes (with hyperbolic attracting cycles) from chaotic
regimes, we conjecture that also the maps in our family\ have a similar
property. As soon as a second fixed point, also virtual, exists for one of the
functions defining the map, we may expect either chaos or some hyperbolic
attracting cycle.

Let us discuss a numerical example. We add an offset $\mu_{R}$ to the function
in the right partition of map $G_{3}$ in (\ref{new}), obtaining the map
\begin{equation}
\widetilde{G}_{3}:\left\{
\begin{array}
[c]{c}%
G_{L}:x^{\prime}=s_{L}x\ \ \ \ \ \ \ \\
G_{M}:x^{\prime}=s_{M}x\ \ \ \ \ \ \ \\
G_{R}:x^{\prime}=s_{R}x+\mu_{R}%
\end{array}
\right.
\begin{array}
[c]{c}%
\text{for} \ x\leq-Z_{L}\\
\text{for} \ -Z_{L}<x<Z_{R}\\
\text{for} \ x\geq Z_{R}%
\end{array}
\label{G3mu}%
\end{equation}
In our financial market model defined in Sec.~2, this map in (\ref{G3mu})
corresponds to a scenario in which the parameter $a^{1,b}$ and/or $a^{2,b}$ of
the trading rule of type 1 chartists differs not only with respect to the sign
from that of type 1 fundamentalists. Let us reconsider the parameters as in
the example shown in Fig.~7. A 2D bifurcation diagram of map 
$\widetilde{G}_{3}$ in (\ref{G3mu}) varying $\mu_{R}$ in the range $(-0.2,0.4)$ and slope
$s_{M}$ in the middle partition in the range $(-0,3,-0.1)$ is shown in
Fig.~11(a). The white region there shown, for $\mu_{R}<0,$ corresponds to the
existence of chaotic dynamics. The colored regions for $\mu_{R}>0$ denote the
existence of attracting cycles. Different colors correspond to cycles of
different periods. The yellow region denotes the existence of an attracting
fixed point.

\begin{figure}
	\begin{center}	
		\includegraphics[width=1\textwidth]{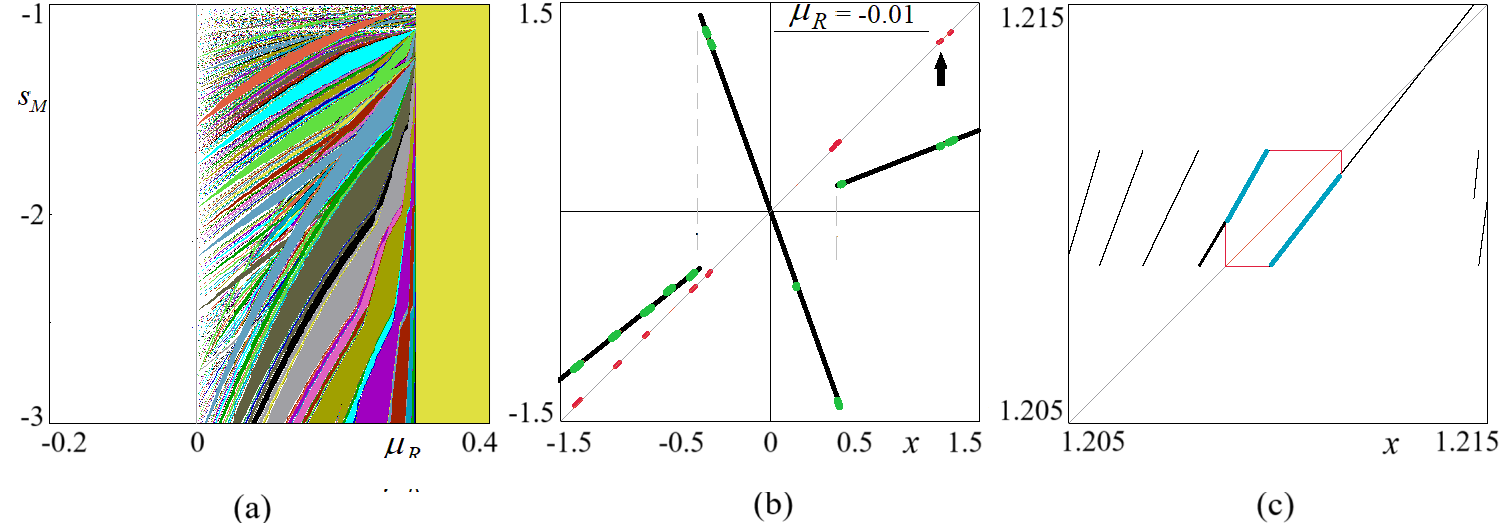} 
		\caption{\label{f11} In (a), 2D bifurcation diagram in the parameter plane ($\mu_{R},s_{M}$)
			of map $\widetilde{G}_{3}$ in (\ref{G3mu}) at $Z_{L}=0.5,$ $Z_{R}=0.5$,
			$s_{L}=0.8,$ $s_{R}=0.4.$ In (b), an example at $\mu_{R}=-0.01$ and
			$s_{M}=-2.82$.\ In (c), the first return map in the interval evidenced by an
			arrow in (b), which is an overlapping PWL Lorenz map, with chaos in the
			interval.}
	\end{center}
\end{figure} 
\begin{figure}
	\begin{center}	
		\includegraphics[width=1\textwidth]{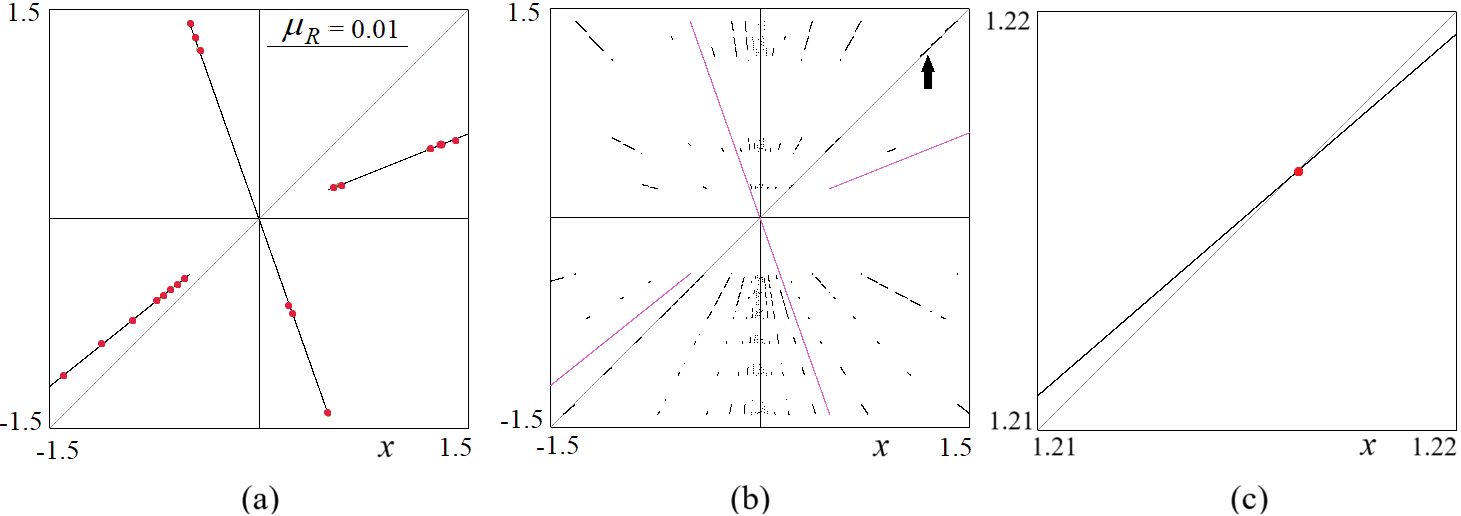} 
		\caption{\label{f12} In (a), map $G_{3}$ in (\ref{new}) with $Z_{L}=0.5,$ $Z_{R}=0.5$,
			$s_{L}=0.8,$ $s_{M}=-2.82,$ $s_{R}=0.4,$ perturbed with $\mu_{R}=0.01$, and
			the trajectories converge to a cycle of period 19. In (b), map $G_{3}$ is
			shown in pink and the 19th iterate of the map is shown in black. It has 19
			segments crossing the main diagonal. An enlargement of the segment evidenced
			by an arrow is shown in (c), where the fixed point $x^{\ast}=1.2161$ is marked
			in red.}
	\end{center}
\end{figure}

An example at $s_{M}=-2.82$ (as in Fig.~7(a,b)) and with $\mu_{R}=-0.01$ is
shown in Fig.~11(b), with dynamics that is chaotic in some intervals, as proven
via the first return map in one of the intervals (evidenced by an arrow in
Fig.~11(b)) and reported in Fig.~11(c), which shows an overlapping Lorenz map,
chaotic in that interval. 

Differently, an example at $s_{M}=-2.82$ (as in Fig.~7(a,b)) and $\mu_{R}=0.01$
is shown in Fig.~12(a), where the dynamics are converging to an attracting
19-cycle. The 19th iterate of the map is shown in Fig.~12(b), evidencing 19
segments crossing the main diagonal; an enlargement is also shown in Fig.~12(c)
with the fixed point in $x^{\ast}=1.2161$. As our simulations suggest, this
cycle is the only stable one, attracting all the points except for the origin
(i.e., there are no other attracting or repelling cycles). 

It may be interesting to compare the dynamics that can be observed in the time
domain, starting from the same initial condition and after the same number of
transient iterations. The three cases with $\mu_{R}=0,$ $\mu_{R}=-0.01$ and
$\mu_{R}=0.01$ are shown in Fig.~13. 

\begin{figure}[ht]
	\begin{center}	
		\includegraphics[width=1\textwidth]{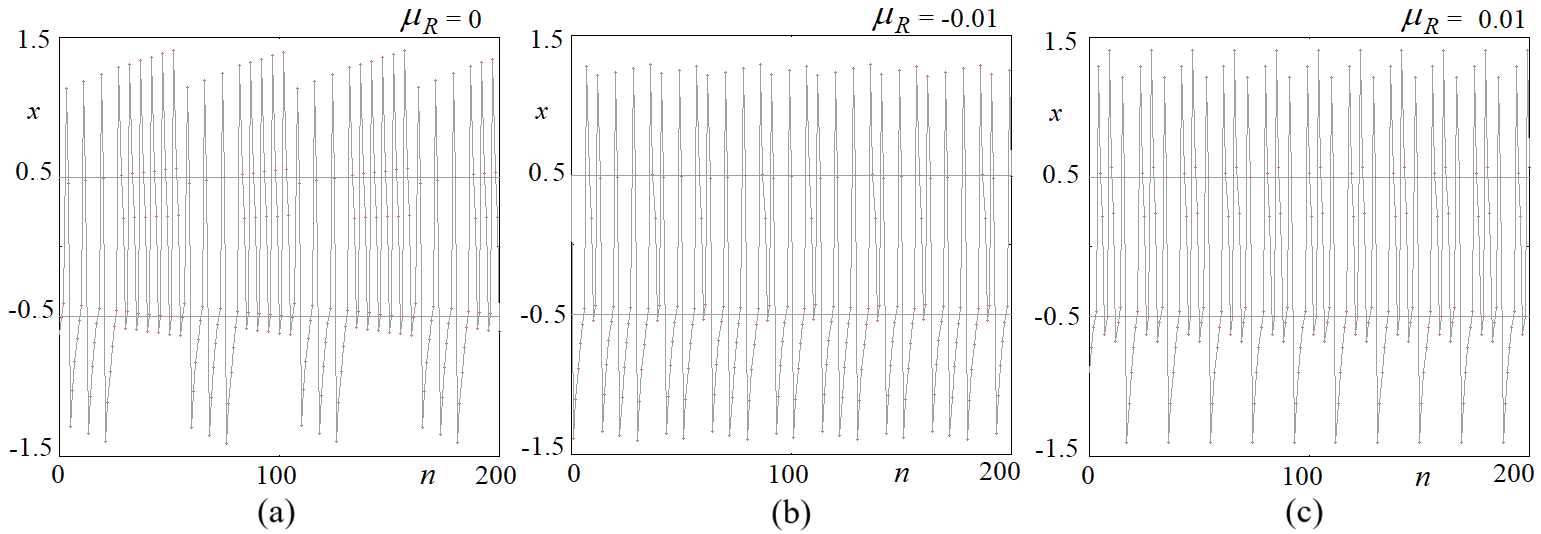} 
		\caption{\label{f13} Map $G_{3}$ in (\ref{new}) $Z_{L}=0.5,$ $Z_{R}=0.5$ and $s_{L}=0.8,$
			$s_{M}=-2.82,\ s_{R}=0.4$ perturbed with $\mu_{R}.$ Versus time trajectory,
			for 200 iterations. In (a) $\mu_{R}=0.$ In (b) $\mu_{R}=-0.01.$ In (c)
			$\mu_{R}=0.01$.}
	\end{center}
\end{figure}

All three parameter constellations can produce endogenous price fluctuations
-- this is a robust feature of our financial market model. At least at first
glance, one may furthermore conclude that the differences in the dynamics are
not so great. However, this is clearly not the case. The dynamics are
quasiperiodic (left), chaotic (middle), or periodic (right). This may have
relevant policy implications, e.g., with respect to stabilization policies.
For instance, certain chaos control methods may help stabilize the dynamics in
the middle panel, but it is not clear whether this will also be the case in
the other parameter constellations. Put differently, it is important to be
able to distinguish between different types of dynamic behaviors. Adding noise
to the dynamics may have further significant implications. In the presence of
a positive Lyapunov exponent, even small shocks may have a marked impact on
the evolution of financial prices, while this is not the case when the
Lyapunov exponent is zero or negative.

\section{Conclusions and future works}

In this work, we have characterized the dynamics that may occur in a wide
class of 1D PWL discontinuous maps with a finite number of discontinuity
points. We have shown that, although numerical simulations may suggest the
existence of chaotic behavior, there is no chaos. The bounded dynamics, not
related to the fixed point, can ultimately be described via a PWL circle map,
so that either nonhyperbolic cycles or quasiperiodic orbits densely fill some
segments. The only attractors different from the fixed point are segments
filled with quasiperiodic trajectories. However, due to the existence of many
partitions, quasiperiodic orbits mimic chaotic behavior, and this, in systems
with practical applications, may still be interpreted as complex behavior.
Moreover, our model, which represents the dynamics of a financial market, may
be used to understand the mechanisms leading to oscillations in prices.

One interesting question left for future work is whether any map in the class
defined here truly acts as a separator between regular and chaotic regimes.
Another interesting question is whether a map satisfying Definition 1 can
appear as a restriction, or first return map, on some segment of a straight
line in a two-dimensional PWL discontinuous map, such as the one proposed in
\cite{GRSSW-24b}. Such systems can be generalized, leading to families of 2D
discontinuous maps whose first return map on some segments of a straight line,
when they exist, lead to a map in our definition. That is, the 1D maps in our
definition occur as particular cases of a larger class of 2D PWL discontinuous
maps. Furthermore, in such cases, it is important to examine whether the role
of these maps as separators between regular and chaotic regimes persists, so
that under perturbation, chaos can emerge, as we observed in the financial
market model when assuming that one\ parameter among $a^{1,b},$ $a^{2,b}$,
$b^{1,d}$, $b^{2,d},$\ in the trading rule of chartists differs from the one
in the trading rule of fundamentalists.

\bigskip

\textbf{Acknowledgements}

Laura Gardini thanks the Czech Science Foundation (Project 22-28882S), the
VSB---Technical University of Ostrava (SGS Research Project SP2024/047), the
European Union (REFRESH Project-Research Excellence for Region Sustainability
and High-Tech Industries of the European Just Transition Fund, Grant
CZ.10.03.01/00/22 003/000004). Davide Radi thanks the Gruppo Nazionale di
Fisica Matematica GNFM-INdAM for financial support. The work of Davide Radi
and Iryna Sushko has been funded by the European Union - Next Generation EU,
Mission 4: "Education and Research" - Component 2: "From research to
business", through PRIN 2022 under the Italian Ministry of University and
Research (MUR). Project: 2022JRY7EF - Qnt4Green - Quantitative Approaches for
Green Bond Market: Risk Assessment, Agency Problems and Policy Incentives -
CUP: J53D23004700008.


\begin{thebibliography}{99}                                                                                               %


\bibitem {Avrutin}Avrutin, V., Gardini, L., Sushko, I., Tramontana, F.:
Continuous and Discontinuous Piecewise-Smooth One-Dimensional Maps. World
Scientific, Singapore (2019)

\bibitem {Berry}Berry, D., Mestel, B.: Wandering interval for Lorenz maps with
bounded nonlinearity, Bull. London Math. Soc., \textbf{23}, 183--189 (1991)

\bibitem {Brianzoni}Brianzoni, S., Campisi, G.: Dynamical analysis of a
financial market with fundamentalists, chartists, and imitators. Chaos,
Solitons \& Fractals, \textbf{130}, 109434,
https://doi.org/10.1016/j.chaos.2019.109434 (2020)

\bibitem {Chiarella-09}Chiarella, C., Dieci, R., He, X.-Z.: Heterogeneity,
market mechanisms, and asset price dynamics. In: Hens, T., Schenk-Hopp\'{e},
K.R. (eds.) Handbook of Financial Markets: Dynamics and Evolution, 277--344.
North-Holland, Amsterdam (2009)

\bibitem {Day}Day, R., Huang, W.: Bulls, bears and market sheep. Journal of
Economic Behavior and Organization, \textbf{14}, 299-329 (1990)

\bibitem {deMelo}de Melo, W., van Strien, S.: One-Dimensional Dynamics,
Springer, New York (1991)

\bibitem {Devaney}Devaney, R.: An introduction to chaotic dynamical systems.
Menlo Park, California: The Benjamin/Cummings Publishing Corporation (1986)

\bibitem {G-2010}Gardini, L., Tramontana, F, Avrutin, V., Schanz, M.: Border
Collision Bifurcations in 1D PWL map and Leonov's approach. International
Journal of Bifurcation \& Chaos, \textbf{20}(10) 3085-3104 (2010)

\bibitem {GT-2012}Gardini, L., Tramontana, F.: Structurally unstable regular
dynamics in 1D piecewise smooth maps, and circle maps. Chaos, Solitons \&
Fractals, \textbf{45}, 1328--1342 (2012)

\bibitem {GRSSW-24}Gardini, L., Radi, D., Schmitt, N., Sushko, I., Westerhoff,
F.: On the limits of informationally efficient stock markets: New insights
from a chartist-fundamentalist model. Submitted for publication.
https://arxiv.org/abs/2410.21198 (2025a)

\bibitem {GRSSW-24b}Gardini, L., Radi, D., Schmitt, N., Sushko, I. and
Westerhoff, F.: On the emergence and properties of weird quasiperiodic
attractors. Submitted for publication. https://arxiv.org/pdf/2503.11264. (2025b)

\bibitem {GRSSW-24c}Gardini, L., Radi, D., Schmitt, N., Sushko, I. and
Westerhoff, F.: How risk aversion may shape the dynamics of stock markets: a
chartist-fundamentalist approach. Working Paper, University of Bamberg. (2025c)

\bibitem {Hommes-09}Hommes, C., Wagener, F.: Complex evolutionary systems in
behavioral finance. In: Hens, T., Schenk-Hopp\'{e}, K.R. (eds.) Handbook of
Financial Markets: Dynamics and Evolution, 217--276. North-Holland, Amsterdam (2009)

\bibitem {HD-93}Huang,W., Day, R.: Chaotically switching bear and bull
markets: the derivation of stock price distributions from behavioral rules.
In: Day, R., Chen, P. (eds.) Nonlinear Dynamics and Evolutionary Economics,
169--182. Oxford University Press, Oxford (1993)

\bibitem {H10}Huang, W., Zheng, H., Chia, W.M.: Financial crisis and
interacting heterogeneous agents. J. Econ. Dyn. Control, \textbf{34},
1105--1122 (2010)

\bibitem {Kollar}Kollar, L.E., Stepan, G., Turi, J.: Dynamics of piecewise
linear discontinuous maps. Int. J. Bif. Chaos, \textbf{14}, 2341-2351 (2004)

\bibitem {Lu}Lu, S., Oberst, S., Zhang, G., Luo, Z.: Bifurcation analysis of
dynamic pricing processes with nonlinear external reference effects. Commun.
Nonlinear Sci. Numer. Simul., \textbf{79}, 104929.S0218127419300271 (2019)

\bibitem {Lux-09}Lux, T.: Stochastic behavioural asset-pricing models and the
stylize facts. In: Hens, T., Schenk-Hopp\'{e}, K.R. (eds.) Handbook of
Financial Markets: Dynamics and Evolution, 161--216. North-Holland, Amsterdam (2009)

\bibitem {Milnor}Milnor, J.: On the concept of attractor. Commun Math Phys,
\textbf{99}, 177--95 (1985)

\bibitem {Rand}Rand, D.:\ The topological classification of Lorenz attractors.
Math. Proc. Camb. Phil. Soc., \textbf{83}, 451--460 (1978)

\bibitem {Segura}Segura, J,. Hilker, F.M., Franco, D.: Degenerate period
adding bifurcation structure of one-diemnsional bimodal piecewise linear maps.
SIAM J. Appl. Math. \textbf{80}(3), 1356-1376 (2020)

\bibitem {TWG-10}Tramontana, F., Westerhoff, F., Gardini, L.: On the
complicated price dynamics of a simple one-dimensional discontinuous financial
market model with heterogeneous interacting traders. J. Econ. Behav. Org.,
\textbf{74}, 187--205 (2010)

\bibitem {TGW-11}Tramontana, F., Westerhoff, F., Gardini, L.: Heterogeneous
speculators and asset price dynamics: further results from a one-dimensional
discontinuous piecewise-linear map. Comput Econ., \textbf{38}, 329--347, DOI
10.1007/s10614-011-9284-9 (2011)

\bibitem {TWG-13}Tramontana, F., Westerhoff, F., Gardini, L.: The bull and
bear market model of Huang and Day: Some extensions and new results. Journal
of Economic Dynamics \& Control, \textbf{37}. 2351--2370 (2013)

\bibitem {TWG-14}Tramontana, F., Westerhoff, F., Gardini, L.: One-dimensional
maps with two discontinuity points and three linear branches: mathematical
lessons for understanding the dynamics of financial markets. Decisions Econ
Finan. \textbf{37, }27--51, DOI 10.1007/s10203-013-0145-y (2014)

\bibitem {TWG-15}Tramontana, F., Westerhoff, F., Gardini, L.: A simple
financial market model with chartists and fundamentalists: Market entry levels
and discontinuities. Math. Comput. Simul. \textbf{108}, 16--40 (2015).

\bibitem {W-09}Westerhoff, F.: Exchange rate dynamics: a nonlinear survey. In:
Rosser Jr, J.B. (ed.) Handbook of Research on Complexity, 287--325. Edward
Elgar, Cheltenham (2009)

\bibitem {Wiggins}Wiggins, S.: Introduction to Applied Nonlinear Dynamical
Systems and Chaos. Springer-Verlag, NY. (1990)

\bibitem {Zhu}Zhusubaliyev, Z.T., Mosekilde, E.: Bifurcations and chaos in
piecewise-smooth dynamical systems. Nonlinear Science A, \textbf{44}, World
Scientific, Singapore (2003)\bigskip
\end{thebibliography}
\end{document}